\documentclass[12pt]{amsart}
\usepackage{amscd,amsthm,amsfonts,amssymb,amsmath,euscript}
\usepackage[dvips]{graphicx}
\usepackage[dvips]{graphics}
\usepackage[T2A]{fontenc}
\usepackage[cp866]{inputenc}

\newtheorem{theorem}{Theorem}[section]
\newtheorem{lemma}{Lemma}[section]
\newtheorem{corollary}{Corollary}[section]
\newtheorem{example}{Example}[section]

\newcommand{\Cov}{\mathop{\sf Cov}\nolimits}
\newcommand{\Aut}{\mathop{\sf Aut}\nolimits}
\newcommand{\tr}{\mathop{\sf tr}\nolimits}
\newcommand{\End}{\mathop{\sf End}\nolimits}

\newcommand{\Rep}{\mathop{\sf Re}\nolimits}
\newcommand{\Imp}{\mathop{\sf Im}\nolimits}

\begin{document}

\title{Hurwitz numbers for regular coverings of surfaces by seamed surfaces and Cardy-Frobenius algebras
of finite groups}

\author{A.V.Alexeevski \and  S.M.Natanzon$^{(*)}$}

\address{A.N.Belozersky Institute,
Moscow State University, Moscow 119899\\}

\email{aba@belozersky.msu.ru}

\address{A.N.Belozersky Institute,
Moscow State University, Moscow 119899\\
Independent University of Moscow\\
Institute Theoretical and Experimental Physics}

\email{natanzon@mccme.ru}
\thanks{$^{(*)}$ Supported
by grants RFBR-07-01-00593, NWO 047.011.2004.026 (RFBR
05-02-89000-HBO-a), NSh-4719.2006.1, INTAS 05-7805, CRDF RM
1-2543-MO-03}

\date{}

\maketitle
          \rightline{Dedicated to S.P.Novicov on the occasion of his 70th birthday}

\begin{abstract}
Analogue of classical Hurwitz numbers is defined in the work for
regular coverings of surfaces with marked points by seamed surfaces.
Class of surfaces includes surfaces of any genus and orientability,
with or without boundaries; coverings may have certain singularities
over the boundary and marked points.  Seamed surfaces introduced
earlier are not actually surfaces. A simple example of seamed
surface is book-like seamed surface: several rectangles glued by
edges like sheets in a book.

We prove that Hurwitz numbers for a class of regular coverings
with action of fixed finite group $G$ on cover space such that
stabilizers of generic points are conjugated to a fixed subgroup $K\subset G$ defines
a new example of Klein Topological Field Theory (KTFT).
It is known that KTFTs are in one-to-one correspondence
with certain class of algebras, called in the work Cardy-Frobenius algebras.
We constructed a wide class of Cardy-Frobenius algebras, including particularly
all Hecke algebras for finite groups.
Cardy-Frobenius algebras corresponding to regular coverings
of surfaces by seamed surfaces are described in terms of group $G$
and its subgroups.
As a result, we give an algebraic formula
for introduced Hurwitz numbers.
\end{abstract}

\tableofcontents

\section{Introduction}

In the paper we construct new examples of Klein Topological Field
Theories (hereinafter abbreviated to KTFT) defined in \cite{AN}. The
examples are connected with Hurwitz numbers for regular coverings of
Klein surfaces by seamed surfaces. These numbers, generalizing
classical Hurwitz numbers, are introduced in the work.

A classical Hurwitz number is a weighted number of non-equivalent
meromorphic maps $f:C\to P$ of Riemann surfaces $C$ to a fixed
Riemann surface $P$  such that topological types of critical values
$p_i\in P$ of $f(z)$ coincide with prescribed data. Initially these
numbers were introduced by Hurwitz \cite{H} for Riemann sphere $P$.
A topological type of a critical value $p_i$ is defined by a
partition $\alpha_i$ of $n$, the degree of $f$, into the unordered
sum of positive integers. Actually, Hurwitz numbers depend on a
Riemann surface $P$ up to isomorphism, a number of critical values
and a set of partitions $\alpha_i$ assigned to points of $P$ that
are critical values of $f$. Clearly, the partitions of $n$ are in
one-to-one correspondence with conjugacy classes of the symmetric
group $S_n$.

Forgetting of the analytic structure, Hurwitz numbers can be defined
as weighted numbers of non-equivalent branched coverings of a fixed
oriented topological surfaces $P$ with prescribed types of branching
at fixed set of marked points. Again, branching type is described by
a conjugacy class of the symmetric group. Moreover, "analytical" and
topological Hurwitz numbers coincide.

Hurwitz numbers can be generalized to the case of regular coverings
of a surface $P$ with marked points \cite{D1,CNP}. Let a finite
group $G$ acts effectively and continuously on a surface $C$. Then
the projection of $C$ onto the orbit space $C/G$ is a branching
covering and the topological type of the branching at a singular
point in $C/G$ is described by a conjugacy class $\alpha$ of $G$.
Let $r:C/G\to P$ be a homeomorphism  such that images of singular
points coincide with marked points of $P$. Then the composite map
$C\to P$ is called a $G$-covering of $P$. Assigning a conjugacy
class $\alpha_i$ to a marked point $p_i$, we can define a Hurwitz
number in the same way as in classical case. If $G=S_n$ then Hurwitz
numbers of $G$-coverings actually coincide with classical Hurwitz
numbers.

In \cite{AN} Hurwitz numbers were extended to a wider category
$\mathcal{K}$ of surfaces. An object of $\mathcal{K}$ is a compact
topological surface of any genus with or without boundary,
orientable or non-orientable and with fixed set of marked points
both in the interior and on the boundary of the surface. A covering
of a surface $\Omega\in\mathcal{K}$ is defined as a continuous map
$f:T\to\Omega$ of a surface $T\in\mathcal{K}$ onto $\Omega$
satisfying the following conditions: $f$ is a local homeomorphism at
interior non-marked points; $f$ is a branching covering in the
neighborhood of an interior marked point; either $f$ is a
homeomorphism or $f$ is topologically equivalent to the complex
conjugation map $z\to \bar z$ in a neighborhood of a point $z\in T$
such that $f(z)$ is boundary non-marked point.

Topological type of a covering $f:T\to\Omega$ at boundary marked
point is defined by an equivalence class $\beta$ of ordered pairs
$(s_1,s_2)$ of two elements of order $1$ or $2$ of symmetric group
$S_n$. The equivalence is defined by $(s_1,s_2)\sim
(gs_1g^{-1},gs_2g^{-1})$ for $g\in S_n$. Actually, Hurwitz numbers
of \cite{AN} depend on a surface $\Omega\in\mathcal{K}$ up to
isomorphism in category $\mathcal{K}$, a set of conjugacy classes
$\alpha_i$ of group $S_n$ assigned to interior marked points of
$\Omega$ and a set $\beta_j$ of equivalence classes of pairs
$(s_1,s_2)$ assigned to boundary marked points.

A real algebraic curve is a pair $(P,\tau)$, where $P$ is a Riemann
surface (representing the complexification of the real algebraic
curve) and $\tau:P\to P$ is an antiholomorphic involution
(representing the complex conjugation on the complexification). The
complex structure on $P$ generates on the factor space
$\Omega=P/\tau$ the structure of a Klein surface \cite{N}.
Forgetting this structure, we obtain that $\Omega$ is a surface from
the category $\mathcal{K}$. Conversely, each $\Omega\in\mathcal{K}$
can be endowed with a structure that comes from a real algebraic
curve. By above, Hurwitz numbers of \cite{AN} can be considered as
Hurwitz numbers for real meromorphic maps with values in a fixed
real algebraic curve.

Next generalization of Hurwitz numbers were done in \cite{AN1,AN2} where coverings of
surfaces $\Omega\in\mathcal{K}$ by so called seamed surfaces were considered.
Combining approaches, in this work we deal with $G$-coverings of surfaces
by seamed surfaces.

Seamed surfaces appear first in physical models \cite{R,KR}. Let
$\Omega_i\in\mathcal{K}$ be a set of surfaces with marked points.
Boundary marked points divide boundary contours of surfaces into
segments. Let a set of homeomorphisms of the segments is given. Then
we can identify a point of a segment with all its images. Obtained
topological space is called a (topological) seamed surface. Images of
$\Omega_i$ in the seamed surface are called regular parts. Note,
that more then two segments of surfaces $\Omega_i$ can be identified
with one segment of the seamed surface. Thus, a seamed surface
generally is not a surface. Several copies of a rectangle glued by
edges like sheets in a book give an example of a seamed surface.

Let finite group $G$ acts effectively and continuously on a seamed surface $T$.
Suppose, $G/T$ is a surface (not seamed surface, but a surfaces of
category $\mathcal{K}$). Then the projection $f:T\to T/G$ is a local
homeomorphism at a generic point $z$ of $T$. Note that the stabilizer $K$ of $z$
is not necessarily trivial subgroup. Clearly, stabilizers of
other generic points of $T$ are conjugated to $K$. Fix one of that
stabilizers, namely, $K$. The topological
type of the covering $f$ at interior singular point $p$ is defined by
a conjugacy class $\alpha$ of the group $N=N_G(K)/K$ where $N_G(K)$ denotes
the normalizer of the subgroup. It is shown in section \ref{local} that
topological type of $f$ at boundary singular point $q$ is described by
an equivalence class $\beta$ of an ordered pair $(S_1,S_2)$ of subgroups
such that $S_1\supset K$, $S_2\supset K$. The equivalence is defined by formula
$(S_1,S_2)\sim (gS_1g^{-1},gS_2g^{-1})$ for $g\in N_G(K)$.

By analogy with the case of coverings of surfaces by surfaces in
category $\mathcal{K}$ we define "topological" Hurwitz numbers for
$G$-coverings of surfaces by seamed surfaces (see section
\ref{cut}).

Following \cite{N1,N2} it is possible to define structures of Klein
surfaces on regular parts of topological seamed surfaces. (The case
of orientable regular parts of a seamed surface and
and Riemann structures on them was considered in \cite{R}.)
Using these structures we
can define "analytical" Hurwitz numbers. It follows from
\cite{N1,N2} that "topological" and "analytical" Hurwitz numbers
coincide. In this work we will explore the topological category of
seamed surfaces.

Topological field theories were introduced in works \cite{At,Se}. In
the case of dimension two it is a system of tensors that are
connected with closed (that is oriented without boundary) surfaces
with marked points. Let $P$ be a such surface and $V$ be a fixed
vector surface. Assign vector space $V_i\approx V$ to each marked
point $p_i$ of $P$ and put $V_P=\otimes_i V_i$. Then the system of
linear maps $\Phi_P:V_P\to\mathbb{C}$ is call a closed topological
field theory if certain axioms are satisfied. It follows from
\cite{D1} that the Hurwitz numbers for the $G$-coverings generate
the tensors that form a closed topological field theory .

Closed  topological field theories are in one-to-one correspondence
with commutative Frobenius algebras equipped with a fixed
nondegenerate symmetric bilinear form and an involutive
antiautomorphism \cite{D2,AN}. For a finite group $G$ this Frobenius
algebra is isomorphic to the center $A$ of the group algebra of $G$
\cite{D1}.

Closed topological field theories were generalized to the case of
oriented surfaces with boundaries \cite{Laz,Moore0}. These theories
are called open-closed topological field theories. Klein topological
field theories (KTFTs)\cite{AN} are generalizations of open-closed
topological field theories to nonorientable surfaces. To define a
KTFT, let us fix two vector spaces, $A$ and $B$. For any
$\Omega\in\mathcal{K}$  we assign a vector space $A_i$ isomorphic to
$A$ to each interior marked point $p_i$ of $\Omega$ and a vector
spaces $B_j$ isomorphic to $B$ to each boundary marked point $q_j$.
KTFT is a system of linear maps (correlators)
$\Phi_\Omega:V_\Omega\to\mathbb{C}$ where $V_\Omega=(\otimes_i
A_i)\otimes(\otimes_j B_j)$, satisfying axioms (see section
\ref{klein_t}). Its restriction to subcategory of orientable
surfaces is open-closed topological field theory.

Relations between correlators and their algebraic implementations
were described in \cite{Laz,Moore0} for open-closed topological
field theory and in \cite{AN} for  KTFT. It was shown in \cite{AN}
that there is a one-to-one correspondence between KTFT and certain
class of algebras. Algebras of that class was called structure
algebras in \cite{AN} and renamed to Cardy-Frobenius algebras in
\cite{AN1,AN2}. Open-closed topological field theories can be
considered as restrictions of KTFTs to the oriented surfaces. The
one-to-one correspondence between open-closed topological field
theories and Frobenius algebras endowed with additional structures
was proved also independently in \cite{Lauda, Moore}.

By definition, a Cardy-Frobenius algebras $H$ is a sum $H=A\oplus B$
of a commutative Frobenius subalgebra $A$ and (typically
noncommutative) Frobenius algebra $B$ which is  a two-sided ideal of
$H$. Algebra $H$ is endowed with additional structures: unit
elements in algebras $A$ and $B$, fixed nondegenerate symmetric
invariant bilinear forms on $A$ and $B$, involutive
antiautomorphisms of $A$ and $B$, element $U\in A$. Main condition
establishes a relation between $A$ and $B$. That relation is an
implementation of Cardy condition (see section \ref{definition_c}).

Semisimple Cardy-Frobenius algebras were classified in \cite{AN}.

In this paper we construct a series of new examples of KTFTs. First,
we prove that Hurwitz numbers for $G$-coverings of surfaces by
seamed surfaces with fixed finite group $G$ and subgroup $K$ which
is the stabilizer of a generic point, generate a Cardy-Frobenius
algebra. Hurwitz numbers of regular coverings in the case $G=S_n$
and $K=\{e\}$ coincides with Hurwitz numbers of coverings of Klien
surfaces by seamed surfaces \cite{AN1,AN2}. The KTFT corresponding
to the latter Hurwitz numbers and Cardy-Frobenius algebra
corresponding to that KTFT were described in \cite{AN1,AN2}. Note,
that noncommutative part $B$ of this Cardy-Frobenius algebra is an
algebra of bipartite graphs.

KTFTs corresponding to $G$-coverings of surfaces by seamed surfaces
are included in the work in broader series of examples of KTFTs: we
prove that a KTFT corresponds to any effective action of a finite
group $N$ on a finite set $X$. Noncommutative parts of
Cardy-Frobenius algebras for those KTFTs include all Hecke algebras
of finite groups. We are grateful to E.Vinberg who pointed to us
this connection. All constructed Cardy-Frobenius algebras $H=A\oplus
B$ are semisimple. To prove it we construct a canonic faithful
representation of  algebra $H$ in vector space $V_X$ of formal
linear combinations of elements of $X$ and prove that the image
coincides with the set of intertwining operators for natural
representation of the group $N$ in $V_X$.

Following \cite{AN}, we obtain finally the formula for
Hurwitz numbers for $G$-coverings of surfaces by seamed surfaces in terms of
Cardy-Frobenius algebra (see section \ref{formula}).

In section \ref{regularcov}  seamed surfaces  and $G$-coverings of
surfaces by seamed surfaces are defined. Local topological
invariants of $G$-coverings are described. Using these invariants
Hurwitz numbers are defined and their topological properties are
described. In section \ref{kleintopo} by means of Hurwitz numbers
for $G$-coverings of surfaces by seamed surfaces we construct a
Klein Topological Field Theory. In section \ref{cardy} we recall a
definition of a Cardy-Frobenius algebra in a bit more invariant form
than in \cite{AN} and associate a Cardy-Frobenius algebra with any
action of a finite group $G$ on a finite set. In section
\ref{formula} we show that the Cardy-Frobenius algebra of
$G$-coverings of surfaces by seamed surfaces is isomorphic to the
Cardy-Frobenius algebra generated by the action of of group
$N=N_G(K)/K$ on a certain finite set. Finally, we give an algebraic
formula for the Hurwitz numbers.

\section{Regular coverings of surfaces by seamed surfaces}
\label{regularcov}
\subsection{Surfaces and seamed surfaces}
{\label{surface}} Two-dimensional compact manifold $\Omega$, not
necessarily orientable or connected, with or without a boundary, and
with finitely many marked points  is called in the work \textit{ a
surface} for short. Thus, the boundary $\partial\Omega$ of a surface
$\Omega$ is either empty or consists of finite number of contours.
Some of the marked points belong to the interior
$\Omega\setminus\partial\Omega$ of $\Omega$ other belong to the
boundary. We require that any boundary contour contains at least one
marked point. An interior not marked point is called simple.
Orientation of a small neighborhood of a point is called \textit{a
local orientation at the point}.

Using marked points we define a stratification of a surface $\Omega$ as follows.
Marked points are considered as $0$-dimensional strata ($0$-strata), $1$-strata are
defined as boundary segments between neighbor boundary marked points at the same
boundary contour,
and $2$-strata are defined as $\Omega^c\setminus(\partial\Omega^c\cup\Omega^c_0)$
where $\Omega^c$ is a connected component of $\Omega$ and  $\Omega^c_0$
denotes the set of marked points in $\Omega^c$.

A homeomorphism $\phi:\Omega'\to\Omega''$ of surfaces is called {\it
an isomorphism} if the image of the set of marked points of
$\Omega'$ coincides with the set of marked points of $\Omega''$.
Therefore, an isomorphism is compatible with the stratifications of
surfaces.

Closed $1$-stratum of a surface $\Omega$  can be
either a segment or a circle, the latter case corresponds to a
boundary contour with just one marked point. Let $\delta_1$ and
$\delta_2$ be two distinct homeomorphic closed $1$-strata of
$\Omega$. By a homeomorphism between $\delta_1$ and $\delta_2$ we
may identify points of $\delta_1$ and $\delta_2$ and  obtain new,
glued topological surface.

We shell consider below stratified topological spaces obtained by
gluing more than two closed $1$-strata
$\delta_1,\delta_2,\delta_3,\dots$  of a surface $\Omega$ into one
stratum by means of homeomorphisms
$\phi_{1,2}:\delta_1\to\delta_2,\phi_{1,3}:\delta_1\to\delta_3,\dots$.
Several such gluings are allowed. The obtained topological space is
called \textit{a seamed surface} (compare with \cite{R}). Seamed
surfaces are $2$-dimensional stratified topological spaces, i.e. all
strata are homeomorphic to topological manifolds of dimension $2$ or
less. Thus, a seamed surface typically is not a surface.

Formal definition of a seamed surface is as follows. First, let us
define its $1$-dimensional analog, that is a topological graph.
\textit{A graph} $\Delta$ is a $1$-dimensional stratified
topological space with finitely many $0$-strata (vertexes) and
finitely many $1$-strata (edges). It is required that the closure of
any edge is homeomorphic either to a closed segment or to a circle,
any vertex belongs to at least one edge, any edge is adjacent to one
or two vertexes. Connectedness of $\Delta$ is not required.
A graph such that  any vertex belongs to
at least two half-edges we call \textit{an s-graph}.
Clearly, the boundary $\partial\Omega$ of a
surface $\Omega$ with marked points is an  s-graph.

A morphisms  of graphs $\varphi:\Delta'\to\Delta''$ is a continuous
epimorphic maps of graphs compatible with the stratification, i.e.
the restriction of $\varphi$ to any open $1$-stratum (interior of an
edge) of $\Delta'$ is a local (therefore, global) homeomorphism with
appropriate open $1$-stratum of $\Delta''$.

Let $(\Omega,\Delta,\varphi)$ be a triple consisting of a surface
$\Omega$ (typically nonconnected), an s-graph $\Delta$ and a
morphism of s-graphs $\varphi:\partial\Omega\to\Delta$. By
identification points of $\partial\Omega$ with their images in
$\Delta$ we obtain a $2$-dimensional stratified topological space
$T$. Clearly, we may and will identify $\Delta$ with the subset of
$T$. $0$-strata of $T$ are vertexes of $\Delta$ and interior marked
points of $\Omega$, its $1$-strata are edges of $\Delta$, its
$2$-strata are $2$-strata of $\Omega$.

A $2$-dimensional stratified topological space $T$
homeomorphic (as a stratified surface) to a topological space obtained
by the described procedure is called \textit{a seamed surface}.
Note, that any surface $\Omega$ (with marked points) is a seamed
surface.

A topological space obtained from a number of copies of rectangles
$R=\{ z\in\mathbb{C} |  |\Rep z|<1, 0\le \Imp z<1\}$ by gluing
segments $(-1,+1)$ of real axis into one segment (like sheets are
glued in a book) we call  \textit{a book-like (open) seamed
surface}.

Let $\Delta$ be a seamed graph. A connected component of
$U\setminus{v}$ where $U$ is appropriately small neighborhood of a
vertex $v\in\Delta$ is contained in one edge; we call it a germ of
the edge. Two germs are called equal if they are adjacent to the
same vertex and have nonempty intersection. Thus, any edge has  two
germs. For seamed surfaces we can similarly define germs of
$2$-strata. We skip the details.

A pair $(v,g(e))$ consisting of a vertex $v$ and a germ $g(e)$  of an
edge adjacent to $v$ is called \textit{a flag}. Note, that if an edge $e$ is
a loop, i.e. both vertices of $e$ coincide with $v$, then there are
\textit{two distinct flags}  corresponding to the pair $(v,e)$. If
an edge $e$ has two distinct vertexes then there is only one flag
corresponding  $(v,e)$.

Let $C(X)$ be a topological cone over
topological space $X$. Then the topological space $C(X)\setminus X$ is called
\textit{an open cone} over $X$.

Let $T$ be a seamed surface,  $\Delta$ be the seamed graph of $T$,
$\Delta_b$ be the set of vertexes of $\Delta$  and $U$ be an
appropriate small neighborhood of a point $z\in T$. Then

\begin{itemize}
\item if $x\in(T\setminus\Delta)$ then $U$ is  homeomorphic to an open  disc;
\item if $x\in\Delta\setminus\Delta_b$ then $U$ is homeomorphic to a
book-like seamed surface;
\item if $x\in\Delta_b$ then $U$ is homeomorphic to an open cone over
a certain connected  graph $\Gamma$; vertexes of $\Gamma$ are in
one-to-one correspondence with flags of $\Delta$ adjacent to $x$ and
edges of $\Gamma$ correspond to germs of $2$-strata adjacent to $x$.
\end{itemize}

\subsection{Coverings and regular coverings}
\label{regular}
Continuous epimorphic map $\varphi:T'\to T''$ of seamed surfaces
$T'$ and $T''$  such that the restriction of $\varphi$ to any open
stratum of $T'$ is a local homeomorphism with a stratum of the same
dimension of $T''$ we call \textit{a covering} of seamed surfaces. A
covering that is a homeomorphism is called \textit{an isomorphism}
of a seamed surface, an isomorphism of a seamed surface with itself
is called \textit{an automorphism}.

Let $G$ be a finite group, $T$ be a seamed surface and
$G\rightarrow\Aut(T)$ be an action of $G$ on  $T$ by automorphisms
of seamed surface $T$. Typically, we assume that the
action is effective. In that case we identify elements of $G$
with their images. The seamed surface $T$ is called a
\textit{$G$-seamed surface} (it would be better to write 'G-(seamed surface)';
nevertheless, we omit the brackets).

An isomorphism $\phi:T'\to T''$ of $G$-seamed surfaces is called
a \textit{$G$-isomorphism} if $g\phi(z)=\phi(gz)$ for any $g\in G$, $z\in T'$.

Let $(T,\Delta,\phi)$ be a $G$-seamed surface. Then, obviously, the
orbit space $T/G$ has the structure of a seamed surface
$(T/G,\Delta/G,\psi)$ where $\psi:\Delta/G\to T/G$ is the induced
map of orbit space of $G$ action on the seamed graph $\Delta$ to the
orbit space $T/G$. Clearly, the map $r:T\to T/G$ is a covering of
seamed surfaces. Below we restrict ourselves to the case of
coverings of surfaces (defined in section \ref{surface}) by seamed
surfaces.

Fix a surface $\Omega$ with marked points and a finite group $G$. A
covering $f:T\to\Omega$ of a surface by $G$-seamed surface is called
\textit{a $G$-covering} if there exists an isomorphism $f_*:T/G\to
\Omega$ such that $f=f_*\circ r$ where $r:T\to T/G$ is the
map of $T$ onto orbit space $T/G$.

Let $f':T'\to \Omega$ and $f'':T''\to \Omega$ be two $G$-coverings
over the surface $\Omega$. An isomorphism $F:T'\to T''$ is called an
\textit{isomorphism of $G$-coverings} if $f'=f''\circ F$ and
$F(g(x))=gF(x)$ for any $x\in T'$ and $g\in G$. Isomorphic
$G$-coverings over $\Omega$ are called \textit{equivalent
coverings}. Denote by $\Aut_Gf$ the group of automorphisms of a
$G$-covering $f$.

Two $G$-coverings $f':T'\to\Omega'$ and $f'':T''\to\Omega''$ are
called \textit{topologically equivalent} if there exists a
$G$-isomorphism of $G$-seamed surfaces $F:T'\to T''$ and an
isomorphism of surfaces $H:\Omega'\to\Omega''$ such that $f''\circ
F=H\circ f'$. Note that if two $G$-coverings are equivalent then
they are topologically equivalent. Converse statement is not
generally true: there exist nonequivalent $G$-coverings $f'$ and
$f''$ of the same surface $\Omega$ which are topologically
equivalent.

Certainly, all introduced notions can be reduced to $1$-dimensional
case, i.e. coverings of $1$-dimensional stratified topological
manifolds - circles, rays and segments - by $1$-dimensional
"seamed" spaces, i.e. seamed graphs. We skip the details.

\subsection{Local topological invariants of $G$-coverings}
\label{local}

Let $f:T\to\Omega$ be a $G$-covering,   $U$ be a small neighborhood of a point $p\in\Omega$.
Denote by $f_U$ the restriction of $f$ to the preimage $\widehat{U}=f^{-1}(U)$.
Clearly, group $G$ acts on $\widehat{U}$ and $\widehat{U}/G \approx U$.
Similar denotations we use for an arbitrary subset $U\subset\Omega$.

Fix a local orientation of $\Omega$ at point $p$. For points of any
type we will describe complete set of invariants of the covering
$f_U$ up to topological equivalence preserving the orientations. We
show below that coverings over a small neighborhood $U$ can be
reduced to coverings over $1$-dimensional spaces: circles, rays and
segments. In the case of one-dimensional coverings topological
equivalence preserving orientations actually coincides with
equivalence of coverings. In this sections topological equivalence
means topological equivalence preserving local orientations.

\smallskip

\textit{Let $p$ be a simple point of $\Omega$}. We may assume that
the boundary $\partial U=\overline{U}\setminus U$ of $U$ is a contour and all
points of $U$ and  $\partial U$ are simple.
In that case the preimage  $\widehat{U}$ of $U$
consists of $n$ discs, $n$ is the degree of the covering. Each of
discs is mapped onto $U$ homeomorphically. Denote by $K$ the stabilizer
of a point $z\in D$ of a disc $D$ that is a connected component of $\widehat{U}$.
Group $K$ acts trivially on $D$ and $\widehat{U}$ is $G$-isomorphic to
the cross-product $G\times_K D$.
By definition, cross-product $G\times_K D$ is a set of equivalence
classes in the direct product $G\times D$; the equivalence is
defined by formula $(g,p)\sim (gk,k^{-1}p)$ for any $k\in K$. The action of group $G$
on $G\times_K U$ is defined by the formula $g(h,p)=(gh,p)$.

$G$-covering  $f_U:\widehat{U}\to U$  is topologically equivalent to
the $G$-covering $G\times_K D\to D$. The stabilizer $K$ depends on
the choice of a disc $D$. Stabilizer $K'$ of point $z'$ from another
disc is conjugated to $K$: $ K'=gKg^{-1}$ where $g\in G$. Thus, the
only local topological invariant of $G$-covering of $\Omega$ at
simple points is a conjugacy class of the stabilizer $K$ of any
preimage of a simple point; $n=|G/K|$ is the degree of the covering.
($|M|$ denotes the cardinality of a set $M$. )
Clearly, stabilizers of all simple points $z\in f^{-1}(\Omega^c)$ where
$\Omega^c$ is a connected component of $\Omega$ are conjugated.

\smallskip

\textit{Let  $p$ be a marked interior point of $\Omega$}. We may
assume that all points of $U$, except $p$, are simple and $\partial
U$ is homeomorphic to a circle. Therefore, $G$-covering
$f_{\partial U}:\widehat{\partial U}\to \partial U$ is topologically
equivalent to a $G$-covering of a circle. Clearly,
$\widehat{\partial U}$ consists of a number of connected components,
each one homeomorphic to a circle. Let $S$ be one of the components.
Fixed local orientation at point $p$ induces the orientation of
$\partial U$ and $S$. Denote by $K$ the stabilizer of a point $z\in
S$. Clearly, stabilizers of all points of $S$ coincide with $K$. A
path around contour $\partial U$ from $f(z)$ to $f(z)$ in the
direction induced by the local orientation at $p$ induces the path
in $S$ from $z$ to a point $\bar z\in f^{-1}(f(z))$ and, therefore,
an element $a\in G$, such that $az=\bar z$. Element $a$ belongs to
the normalizer $N_G(K)$ and is defined up to the stabilizer $K$.
Denote by $K_a$ the subgroup of $G$ generated by $K$ and $a$. Group
$K_a$ acts on $S$ with kernel $K$,  and, clearly, $S/K_a\approx
\partial U$. Therefore, $G$-covering $f_{\partial
U}:\widehat{\partial U}\to \partial U$ is topologically equivalent
to the $G$-covering $G\times_{K_a} S\to\partial U$.

If we take another connected component $S'$ of $\widehat{\partial
U}$ and another point $z'\in S'$ then we obtain another pair
$(K',a')$ where $K'$ is the stabilizer of $z'$ and $a'$ is element
of the normalizer $N_G(K')$. Clearly, there exists an element $g\in
G$ such that $K'=gKg^{-1}$ and $a'=gag^{-1}$. Thus, the only
topological invariant of $G$-covering over the boundary of $U$ is a
conjugacy class $\alpha$  of a pair $(K,a)$ where $K$ is a
subgroup of $G$ and $a$ is an element of group $N=N_H(K)/K$.

Note that classes of topological equivalence of $G$-covering
$f_{\partial U}$ of a boundary $\partial U$ are in one-to-one
correspondence with classes of topological equivalence of
$G$-coverings $f_U$ of $U$ because $U$ and $\widehat U$ can be
reconstructed from $\partial U$ and $\widehat{\partial U}$ by
canonical operations: $U$ is homeomorphic to an open cone over
$\partial U$ and any connected component of $\widehat{U}$ is
homeomorphic to an open cone over a connected component of
$\widehat{\partial U}$.

Denote by $\mathcal{A}$ the set of conjugacy classes of elements
of group $N=N_G(K)/K$. It was shown that the topological invariant of a $G$-seamed
covering at an interior marked point is an element of  $\mathcal{A}$.
An element of $\mathcal{A}$ we call 'an interior field' and denote usually by $\alpha$.

Define an involute map $\mathcal{A}\to\mathcal{A}^*$ by formula
$[a]\to [a]^*=[a^{-1}]$ where $[a]$ denotes the conjugacy class of
element $a\in N$.

Elements of the centralizer $C_N(a)$ are in one-to-one correspondence
with $G$-automorphisms of the $G$-covering $f_U:\widehat{U}\to U$.
By this reason we define automorphism group of $\alpha\in\mathcal{A}$
as $\Aut\alpha=C_N(a)$ for an element $a\in\mathcal{A}$. The group
$\Aut\alpha$ is defined up to isomorphism.

\smallskip

\textit{Let $p$ belongs to an open $1$-stratum of $\Omega$}. We may
assume that small neighborhood $U$ of $p$ is homeomorphic to a
rectangle $R=\{ z\in\mathbb{C} |  |\Rep z|<1, 0\le \Imp z<1\}$,
point $p$ corresponds to $0$ and $U\cap\partial\Omega$ corresponds
to the segment $(-1,+1)$ of real axis.  The preimage $\widehat{U}$
of $U$ consists of a number of book-like seamed surfaces
$B^1,B^2,\dots$. Denote by $R$ a ray $R\subset U$ that corresponds
to the segment $[0,i)$ of imaginary axis. Then $G$-covering
$f_R:\widehat{R}\to R$ is topologically equivalent to the
$G$-covering of a ray. Choose a ray $\delta\subset\widehat{R}$.
Denote by $K$ common stabilizer of interior points of $\delta$ and
by $S$ the stabilizer of the vertex $z$ of $\delta$. Clearly,
$K\subset S$. If we choose another ray instead of $\delta$ than we
obtain conjugated pair $(K', S')$ of subgroups: $K'=gKg^{-1}$,
$S'=gS g^{-1}$ for some $g\in G$. Up to a topological equivalence
$G$-covering of a ray is defined by conjugacy class of pair of
subgroups $(K,S)$. Indeed, group $S$ acts on the connected component
$\widehat{R}^c$ of $\widehat{R}$ that contains $\delta$ and
$\widehat{R}^c/S\approx R$. $\widehat{R}^c$ consists of $|S/K|$ rays
with common vertex. Therefore, $G$-covering $f_R:\widehat{R}\to R$
is topologically equivalent to the $G$-covering $G\times_S
\widehat{R}^c$.

Classes of topological equivalence of $G$-covering $f_R$ of ray $R$
are in one-to-one correspondence with classes of topological
equivalence of $G$-coverings $f_U$ of $U$ because $U$ and
$\widehat{U}$ can be reconstructed from $R$ and $\widehat{R}$ by
canonical operations: $U$ is homeomorphic to an open cylinder over
$R$ and any connected component of $\widehat{U}$ is homeomorphic to
an open cylinder over a connected component of $\widehat{R}$.

Thus, the only topological invariant at point $p$ is a pair  of
subgroups $(K,S)$ such that $K\subset S$ up to conjugation by
element of $G$. Fix a subgroup $K$ from the pair. Then the set of
conjugacy classes of pairs is in one-to-one correspondence with
classes of equivalence of subgroups $S\subset G$ such that $S\supset
K$; the equivalence is defined by $S\sim nSn^{-1}$ where $n\in
N_G(K)$.

\smallskip

\textit{Let $p$ be a boundary marked point}. We may assume that the boundary
 $\partial{U}=\widehat{U}\setminus U$ of appropriate small neighborhood of $p$
is homeomorphic to a closed segment and  $U$ is homeomorphic to an
open cone over $\partial{U}$. Restriction $f_{\partial{U}}$ of
$G$-covering $f$ to $\partial{U}$ is a $G$-covering of a segment. A
connected component $\widehat{\partial{U}}^c$ of the preimage
$\widehat{\partial{U}}$ is a connected bipartite graph. Choose an
edge $e\in \widehat{\partial{U}}^c$ and denote vertexes of $e$ by
$a$ and $b$. Vertexes are ordered according to the local orientation
at $p$. Let $K$ be common stabilizer of internal points of $e$,
$S_a$ be the stabilizer of $a$ and $S_b$ be the stabilizer of $b$.
Clearly, $K\subset S_a\cap S_b$. Denote by $S$ a subgroup generated
by subgroups $S_a$ and $S_b$. It can be shown that $S$ acts on the
component $\widehat{\partial{U}}^c$ and $\widehat{\partial{U}}^c/S
\approx \partial{U}$. Moreover, the $G$-covering
$f_{\partial{U}}:\widehat{\partial{U}}\to \partial{U}$ is
topologically equivalent to the $G$-covering of cross-product
$G\times_S \widehat{\partial{U}}^c$. Covering over $\partial{U}$
defines the covering $f_U$ of the neighborhood $U$ up to topological
equivalence because $U$ is homeomorphic to an open cone over
$\partial{U}$ and $\widehat{U}$ is homeomorphic to the disjoint
union of open cones over connected components of
$\widehat{\partial{U}}$.

Thus, up to topological invariance, $G$-covering in the neighborhood
of a marked boundary point is defined by an ordered triple of
subgroups $(S_a, K,S_b)$ such that $K\subset S_a\cap S_b$ up to
equivalence $(S_a, K,S_b)\sim (gS_a g^{-1},gKg^{-1},gS_bg^{-1})$ for
$g\in G$. Fix subgroup $K$ from a triple. Clearly, the set
of equivalent classes of triples $(S_a, K,S_b)$ are in one-to-one
correspondence with equivalent classes of pairs $(S_a,S_b)$ such
that $S_a\supset K$ and $S_b\supset K$; the equivalence is defined
by $(S_a, K,S_b)\sim (nS_a,n^{-1},nS_bn^{-1})$, $n\in N_G(K)$. We denote the latter
set by $\mathcal{B}$. An element of $\mathcal{B}$ we call 'a
boundary field' and denote usually by $\beta$.

We may assume that the group $N=N_G(K)/K$ acts on the set of
subgroups $S$ containing $K$ by conjugations $S\to nSn^{-1}$, $n\in
N$ because $K$ acts trivially. Elements of $N$ such that
$nS_an^{-1}=S_a$, $nS_bn^{-1}=S_b$ are in one-to-one correspondence
with $G$-automorphisms of the $G$-covering $f_U:\widehat{U}\to U$.
By this reason we define automorphism group of $\beta\in\mathcal{B}$
as $\Aut\beta=\{n\in N| nS_an^{-1}=S_a, nS_bn^{-1}=S_b\}$ for a pair
of subgroups $(S_a,S_b)\in\beta$. The group $\Aut\beta$ is defined
up to isomorphism.

Define an involute map $\mathcal{B}\to\mathcal{B}$ by formula
$[(S_a,S_b)]\to [(S_a,S_b)]^*=[(S_b,S_a)]$ where $[(S_a,S_b)]$
denotes the equivalence class of an ordered pair of subgroups.

\smallskip

Note that if $\Omega$ is connected surface and $f:T\to\Omega$ is a
$G$-covering then stabilizers of all simple points of $T$ are
conjugated subgroups. In the case of nonconnected $\Omega$  we
restrict ourselves to $G$-coverings such that stabilizers of simple
points over all components of $\Omega$ are conjugated.

\subsection{Hurwitz numbers of $G$-coverings}
\label{cut}
Fix a surface $\Omega$, a finite group $G$ and
a subgroup $K$ of $G$.
Denote by $\Omega_a$ the set of interior marked points of $\Omega$ and by $\Omega_b$ the set
of boundary marked points.

Let $f:T\to \Omega$ be a $G$-covering of $\Omega$ by $G$-seamed
surface such that stabilizers of all simple points of $T$ are
conjugated to $K$. It was shown in section \ref{local} that local
topological invariants of $f$ at interior marked points belong to
the set $\mathcal{A}$ of conjugacy classes of the group $N_G(K)/K$
and  local topological invariants of $f$ at boundary marked points
belong to the set $\mathcal{B}$ consisting of equivalence classes of
ordered pairs $(S_a,S_b)$ of subgroups such that $K\subset S_a\cap
S_b$; the equivalence is defined as conjugation of a pair by an element
$n\in N_G(K)$. Therefore, $G$-covering $f:T\to\Omega$ defines maps
$\mathfrak{a}_f:\Omega_a\to\mathcal{A}$ and
$\mathfrak{b}_f:\Omega_b\to\mathcal{B}$, the images
$\mathfrak{a}_f(p)$ (resp., $\mathfrak{b}_f(q)$) we denote by
$\alpha_p$ (resp.$\beta_q$) and call interior (resp., boundary)
fields.

Fix arbitrary  maps $\mathfrak{a}:\Omega_a\to\mathcal{A}$ and
$\mathfrak{b}:\Omega_b\to\mathcal{B}$. Denote by
$\Cov(\Omega,\mathfrak{a},\mathfrak{b})$ the set of isomorphic
classes of $G$-coverings of $\Omega$ by $G$-seamed surfaces  such
that its local invariant at any point $p\in\Omega_a$ is equal to
$\mathfrak{a}(p)$ and its local invariants at any point
$q\in\Omega_b$ is equal to $\mathfrak{b}(q)$.
\smallskip

The number $\mathcal{H}(\Omega;\mathfrak{a};\mathfrak{b})=
\sum_{[f]\in\Cov(\Omega,\mathfrak{a},\mathfrak{b})}\frac{1}{|\Aut_G
f|}$ where $f$ is a representative of the isomorphic class $[f]$ of
$G$-coverings is called {\it a Hurwitz number of $G$-coverings}.
\smallskip

There are relations between classical Hurwitz numbers for different
Klein surfaces. Those relations correspond to
cuts of surfaces. Similar relations we establish below
for Hurwitz numbers of $G$-coverings by seamed surfaces.

Let $\Omega$ be a surface. A
connected closed oriented curve $\gamma\subset\Omega$ without
self-intersections is called \textit{a simple cut} if it does not
meet any marked point and is  either a closed contour belonging to
the interior of the surface or a segment such that
$\partial\gamma=\gamma\cap\partial\Omega$.

Let $\gamma\subset\Omega$ be a simple cut. We will define a cut
contracted surface $\Omega_\#=\Omega/\gamma$ in two steps. First, a
cut surface $\Omega_*$ is a surface
$(\Omega\setminus\gamma)\cup\widehat{\gamma}$ where
$\widehat{\gamma}$ is the set of pairs $(r,c)$, $r\in\gamma$, $c$ is
a co-orientation of $\gamma$ at point $r$. It is clear how to define
the topology on $\Omega_*$, see \cite{AN} for the details. Clearly,
there is a natural continuous map $\Omega_*\to\Omega$, its
restriction to $\Omega_*\setminus\widehat{\gamma}$ is the
homeomorphism with $\Omega\setminus\gamma$, its restriction to
$\widehat{\gamma}$ is a double-covering of $\gamma$.
Second, {\it a cut contracted surface $\Omega_\#$} is
a surface obtained by contracting each connected component of
$\widehat{\gamma}$ into a point. Here $\widehat{\gamma}$ is
the preimage of simple cut $\gamma$ in cut surface $\Omega_*$.
 Marked points of $\Omega_\#$ are
marked points coming from $\Omega$ and points that are
contracted components of $\widehat{\gamma}$.
Note, that the orientation of $\gamma$ induces local orientations at
corresponding marked points.

Clearly, $\Omega_\#$ canonically corresponds to $\Omega$ despite of the absence
of natural continuous map of $\Omega$ to $\Omega_\#$ (the homeomorphism
between $\Omega\setminus\gamma$ and $\Omega_\#\setminus\mbox{(image of $\widehat{\gamma}$)}$
cannot be extended to $\gamma$).

Depending on a simple cut $\gamma$, three situations may occur.
\\(i)  $\gamma$ is co-orientable contour. Then $\widehat{\gamma}$
consists of two contours and $(\Omega_\#)_a=\Omega_a\sqcup \{p',p''\}$
where points $p^{(i)}$ correspond to components of $\widehat{\gamma}$;
$(\Omega_\#)_b=\Omega_b$.
\\(ii) $\gamma$ is non-coorientable contour. Then $\widehat{\gamma}$
consists of one contour and $(\Omega_\#)_a=\Omega_a\sqcup \{p'\}$
where $p'$ corresponds to $\widehat{\gamma}$; $(\Omega_\#)_b=\Omega_b$.
\\(iii) $\gamma$ is a segment. Then $\widehat{\gamma}$
consists of two segments and $(\Omega_\#)_b=\Omega_b\sqcup \{q',q''\}$
where points $q^{(i)}$ correspond to components of $\widehat{\gamma}$;
$(\Omega_\#)_a=\Omega_a$.

Marked points $p',p''$ (in the case (i)), $p$ (in the case (ii)) and
$q',q''$ (in the case (iii)) we call 'extra marked points'.

Let
$\mathfrak{a}:\Omega_a\to\mathcal{A}$ and
$\mathfrak{b}:\Omega_b\to\mathcal{B}$ be two maps of sets of marked
points of a surface $\Omega$. We will extend them to obtain maps
$\mathfrak{a}_{\Omega_\#}:(\Omega_\#)_a\to\mathcal{A}$ and
$\mathfrak{b}_{\Omega_\#}:(\Omega_\#)_b\to\mathcal{B}$ for cut
contracted surface $\Omega_\#$. The extended maps depend on an
element $\alpha\in\mathcal{A}$ for simple cuts of types (i) and (ii)
and on an element $\beta\in\mathcal{B}$ for a simple cut of type
(iii). For marked points of $\Omega_\#$  coming from the marked
points of $\Omega$ extended maps coincide with maps $\mathfrak{a}$
and $\mathfrak{b}$. Values of maps $\mathfrak{a}_{\Omega_\#}$ and
$\mathfrak{b}_{\Omega_\#}$ at all extra marked points are equal to
$\alpha$  (in cases (i), (ii)) and $\beta$ (in case (iii)).

Identifying maps with sets of their images, we use below the following denotations.
\\(i)\ \  $\mathfrak{a}_{\Omega_\#}=(\mathfrak{a},\alpha,\alpha)$
\\(ii)\  $\mathfrak{a}_{\Omega_\#}=(\mathfrak{a},\alpha)$
\\(iii) $\mathfrak{b}_{\Omega_\#}=(\mathfrak{b},\beta,\beta)$

The relations between Hurwitz numbers are described in the following important lemma.

\begin{lemma} \label{lemma_cut}
Let $\gamma\subset\Omega$ be a simple cut of a surface
$\Omega$ and $\Omega_\#=\Omega/\gamma$  be cut contracted surface.

(i) If $\gamma$ is a co-orientable contour, then

$$\mathcal{H}(\Omega;\mathfrak{a};\mathfrak{b})=
\sum_{\alpha\in\mathcal{A}}|\Aut\alpha|
\mathcal{H}(\Omega_\#;\mathfrak{a},\alpha,\alpha;\mathfrak{b})$$
where group $\Aut\alpha$ is defined in section \ref{local}.

(ii) If $\gamma$ is a non-coorientable contour, then

$$\mathcal{H}(\Omega;\mathfrak{a};\mathfrak{b})=
\sum_{\alpha\in\mathcal{A}} d^\alpha
\mathcal{H}(\Omega_\#;\mathfrak{a},\alpha;\mathfrak{b})$$ where
$d^\alpha$ is the number of elements $n\in N=N_G(K)/K$ such that
$n^2=a^{-1}$ for a fixed representative $a\in\alpha$ of the
conjugacy class $\alpha$\footnote{We correct here an incorrectness
in the definition of numbers $d^\alpha$ in \cite{AN}}.

(iii) If $\gamma$ is a segment, then

$$\mathcal{H}(\Omega;\mathfrak{a};\mathfrak{b})=
\sum_{\beta\in\mathcal{B}}|\Aut\beta|
\mathcal{H}(\Omega_\#;\mathfrak{a};\mathfrak{b},\beta,\beta)$$
where group $\Aut\beta$ is defined in section \ref{local}.

\end{lemma}

Actually, the proof follows the proof of similar lemma
for coverings of surfaces by surfaces \cite{AN}. $\Box$

\section{Klein Topological Field Theory for $G$-coverings
of surfaces by seamed surfaces}
\label{kleintopo}
We reproduce below the definition of Klein topological field theory
given in \cite{AN}.

\subsection{Category of surfaces and Klein functor}
\label{klein_f}

In this section we define a category $\mathcal S$ of surfaces. An
object of $\mathcal S$ is a surfaces $\Omega$  with marked points
(see section \ref{surface}) endowed with a set of local orientations
at marked points. The local orientation at marked point
$r\in\Omega_0$ is denoted by $o_r$.

We define morphisms of four types. First, an isomorphism of surfaces
(see section \ref{surface}) preserving local orientations at marked
points  is called a morphism in the category $\mathcal S$. Second,
for an oriented cut $\gamma\subset\Omega$ the  correspondence
$\Omega\to\Omega_\#$ where $\Omega_\#=\Omega/\gamma$ is the cut
contracted surface is
called a morphism. Third, changing a local orientation at a marked
point $r$ to the opposite orientation is called a morphism. Fourth,
marking off a new point and fixing a local orientation at it is
called a morphism. An arbitrary morphism in $\mathcal S$ is defined
as a combination of morphisms of those four types. The disjoint
union $\Omega'\cup \Omega''$ of two surfaces we consider as a tensor
product in $\mathcal{S}$; all axioms of tensor category are
satisfied.

Let $\{X_m | m\in M\}$ be a finite set of $n=|M|$ vector spaces
$X_m$ over the field of complex numbers $\mathbb{C}$. The action of
the symmetric group $S_n$ on $\{1,\dots,n\}$ induces its action on
the sum of the vector spaces $\left(\oplus_{\sigma}
X_{\sigma(1)}\otimes\dots\otimes X_{\sigma(n)}\right)$ where
$\sigma$ runs through the bijections $\{1,\dots, n\}\to M$; an
element $s\in S_n$ takes $X_{\sigma(1)}\otimes\dots\otimes
X_{\sigma(n)}$ to $X_{\sigma(s(1))}\otimes\dots\otimes
X_{\sigma(s(n))}$. Denote by $\otimes_{m\in M} X_m$ the subspace of
all invariants of this action. The vector space $\otimes_{m\in M}
X_m$ is canonical isomorphic to tensor product of all $X_m$ in any
fixed order; the isomorphism is the projection of $\otimes_{m\in M}
X_m$ to the summand that is equal to the tensor product of $X_m$ in
that order.

Assume that all $X_m$ are equal to a fixed vector space X. Then any
bijection $M \to M'$ of sets induces an isomorphism $\otimes_{m\in
M} X\leftrightarrow \otimes_{m'\in M'} X$.

Fix a pair $(A, B)$ of vector spaces over $\mathbb{C}$, a set of
elements $1_A\in A$, $1_B\in B$, $K_A^\otimes\in A\otimes A$,
$K_B^\otimes\in B\otimes B$, $U\in A$ and involute linear maps
$\iota_A:A\to A$, $\iota_B:B\to B$ such that $\iota_A(1_A)=1_A$,
$\iota_B(1_B)=1_B$. Denote by $\mathcal{N}$ the set of all those
data:  $\mathcal{N}= \{A, B,1_A\in A, 1_B\in B, K_A^\otimes\in
A\otimes A, K_B^\otimes\in B\otimes B, U\in A, \iota_A,\iota_B\}$.
Using $\mathcal{N}$,  we will define a functor $\mathcal{V}$ from
the category $\mathcal{S}$ to the category of vector spaces.

Let $\Omega$ be an object of $\mathcal{S}$, $\Omega_0$ be the set of marked points,
$\Omega_a\subset\Omega_0$ be the set of interior marked points and
$\Omega_b=\Omega_0\setminus\Omega_a$ be the set of boundary marked points.
Assign a copy $A_p$ of vector space $A$ to any $p\in\Omega_a$ and
a copy $B_q$ of vector space $B$ to any $q\in\Omega_b$.
Denote by $V_\Omega$ vector space
$(\otimes_{p\in\Omega_a}A_p)\otimes(\otimes_{q\in\Omega_b}B_q)$
By definition, put $\mathcal{V}(\Omega)=V_\Omega$

Functor $\mathcal{V}$ takes morphisms of four types in the category $\mathcal{S}$
to the following morphisms of linear spaces. For any morphism $\phi$
in category $\mathcal{S}$ the morphism $\mathcal{V}(\phi)$ is denoted below by $\phi_*$.

\smallskip

1) Let $\phi:\Omega'\to\Omega''$ be an isomorphism of surfaces. Then
$\phi_*: V_{\Omega'}\to V_{\Omega''}$ is a linear map induced by the bijections
$\phi|_{\Omega'_a}:\Omega'_a\to\Omega''_a$ and
$\phi|_{\Omega'_b}:\Omega'_b\to\Omega''_b$.

\smallskip

2) Let $\psi:\Omega'\to\Omega''$ be a morphism of changing the local
orientation at an interior point $r\in\Omega_0$. Then
$\psi_*$ is a liner map that is identical on all
components of the tensor product
$V_{\Omega'}=(\otimes_{p\in\Omega_a}A_p)\otimes(\otimes_{q\in\Omega_b}B_q)$
except $A_r$ and it coincides with $\iota_A$ on $A_r$.

The morphism $\psi_*$ for morphism  $\psi$ of changing the local orientation
at a boundary marked point is defined similarly, just replace $\iota_A$ by $\iota_B$.

\smallskip

3) Let $\xi:\Omega'\to\Omega''$ be a morphism of marking off an
interior point $p\in\Omega'\setminus\Omega'_0$ and fixing a local
orientation $o_p$. Vector space $V_{\Omega''}$ can be canonically
identified with $V_{\Omega'}\otimes A$. Then
$\xi_*: V_{\Omega'}\to V_{\Omega''}$ maps vector
$x\in V_{\Omega'}$ to the vector of $V_{\Omega''}$ that corresponds
to $x\otimes 1_A\in V_{\Omega'}\otimes A$. Property
$\iota_A(1_A)=1_A$ guaranties the correctness.

Similarly, the morphism $\xi'_*$ corresponding to a
morphism $\xi'$ which is marking off a boundary
point $q\in\partial\Omega'\setminus\Omega'_b$
takes vector $x\in V_{\Omega'}$ to the vector of $V_{\Omega''}$ corresponding to
$(x\otimes 1_B)\in V_{\Omega'}\otimes B$.

\smallskip

4) Let $\eta:\Omega\to\Omega_\#$ be a morphism from $\Omega$ to cut
contracted surface $\Omega_\#$ induced by an oriented cut
$\gamma\subset\Omega$.  As it was shown in section \ref{cut}, $\eta$ induces
the inclusion map of the set $\Omega_0=\Omega_a\cup\Omega_b$ of
marked points of $\Omega$ into the set of marked points
$\Omega_{\#0}$ of $\Omega_\#$ and
$\Delta=\Omega_{\#0}\setminus\Omega_0$ consists of either one point
or two points depending on $\gamma$.

If $\gamma$ is non-coorientable contour, then  $\Delta$ consist of
one point, $V_{\Omega_\#}$ is canonically isomorphic to
$V_\Omega\otimes A$ and $\eta_*(x)=x\otimes U$. If
$\gamma$ is a co-orientable contour then $|\Delta|=2$,
$V_{\Omega_\#}\equiv V_\Omega\otimes A\otimes A$ and
$\eta_*(x)=x\otimes K_A^\otimes$. If $\gamma$ is a
segment then $|\Delta|=2$, $V_{\Omega_\#}\equiv V_\Omega\otimes
B\otimes B$ and $\eta_*(x)=x\otimes K_B^\otimes$.

\smallskip

In addition, the tensor product in $\mathcal{S}$ defined by the
disjoint union of surfaces $\Omega'\otimes\Omega''\to
\Omega'\cup\Omega''$ induces the tensor product of vector spaces
$\theta_*=V_{\Omega'}\otimes V_{\Omega''}\to
V_{\Omega'\cup\Omega''}$.

\smallskip

The functorial properties of $\mathcal{V}$ can be easily verified.

\subsection{Klein Topological Field Theory}
\label{klein_t}
Fix a set of data
$\mathcal{N}=\{A, B, 1_A\in A, 1_B\in B, K_A^\otimes\in A\otimes A, K_B^\otimes\in B\otimes B, U\in A,
\iota_A:A\to A, \iota_B:B\to B\}$
that defines a Klein functor $\mathcal{V}$.

A system of linear forms $\mathcal {F}= \{\Phi_{\Omega}:V_\Omega\to
\mathbb{C}\}$, $\Omega\in\mathcal{S}$, is
called a \textit{Klein topological field theory (KTFT)} if it
satisfies the following axioms.

$1^\circ$ {\it Topological invariance.}

$$\Phi_{\Omega}(\phi_*(x))=\Phi_\Omega(x)$$

for any isomorphism $\phi:\Omega\to \Omega'$ of surfaces.

\vskip 0.6cm

$2^\circ$ {\it Invariance under a change of local orientations.}

$$\Phi_{\Omega'}(\psi_*(x))=\Phi_\Omega(x)$$ for any morphism
$\psi:\Omega\to\Omega'$ of changing the local orientation at a
marked point.

$3^\circ$ {\it Nondegeneracy.}

Define bilinear form $(x',x'')_A$ on vector space $A$ by formula
$(x',x'')_A=\Phi_{(S,p_1,p_2)}(x'_{p_1}\otimes x''_{p_2})$ where
$(S,p',p'')$ is a sphere with two marked points $p'$ and $p''$,
local orientations at $p',p''$ are induced by a global orientation
of the sphere, $x'_{p_1}\otimes x''_{p_2}$ is the element of
$V_{(S,p_1,p_2)}$ corresponding to $x'\otimes x''\in A\otimes A$
under the canonical isomorphism $V_{(S,p_1,p_2)}\equiv A\otimes A$.
It follows from axioms $1^\circ$ and $2^\circ$ that the definition
is correct and bilinear form $(x',x'')$ is symmetric. Taking a disc
with two boundary marked points $(D,q_1,q_2)$ instead of
$(S,p_1,p_2)$ we define a symmetric bilinear form  $(x',x'')_B$ on
vector space $B$ by similar construction.

The axiom is:

\begin{center}both forms $(x',x'')_A$ and  $(x',x'')_B$ are nondegenerate.\end{center}

$4^\circ$ {\it Invariance under addition of marked point.}

$$\Phi_{\Omega'}(\xi_*(x))=\Phi_\Omega(x)$$

for any morphism $\xi:\Omega\to\Omega'$ that is induced by
marking off an additional point in $\Omega$.

\vskip 0.6cm

$5^\circ$ {\it Cut invariance.}

$$\Phi_{\Omega_\#}(\eta_*(x))=\Phi_{\Omega}(x)$$

for any morphism $\eta:\Omega\to\Omega_\#$ of $\Omega$ to
the cut contracted surface  $\Omega_\#$ that is induced
by an oriented simple cut $\gamma\subset\Omega$.

$6^\circ$ {\it Multiplicativity.}

$$\Phi_{\Omega}(\theta_*(x'\otimes x'))= \Phi_{\Omega'}(x')\Phi_{\Omega''}(x'')$$

for $\Omega=\Omega'\cup\Omega''$,  $x'\in V_{\Omega'}$, $x''\in V_{\Omega''}$
and isomorphism  $\theta^*: V_{\Omega'}\otimes V_{\Omega''}\to V_{\Omega}$
corresponding to the
tensor product  $\theta: \Omega'\times\Omega''\to \Omega$.

\vskip 0.6cm

\begin{lemma} Let $\{\Phi_{\Omega}\}$ be a KTFT. Then the following
relations hold:
\\(i)  $(1_A,\alpha)_A=\Phi_{S^2,p}(\alpha)$ where $S^2$ is a sphere with marked
point $p$, a local orientation at $p$ is induced by a global
orientation of $S^2$, and $\alpha\in A$ is an arbitrary element;
\\(ii) $(1_B,\beta)_B=\Phi_{D^2,q}(\beta)$ where $D^2$ is a disc with one boundary marked
point $q$, a local orientation at $q$ is induced by a global
orientation of $D^2$ and $\beta\in B$ is an arbitrary element;
\\(iii) $(U,\alpha)_A=\Phi_{P^2,p}(\alpha)$ where $P^2$ is a projective plane with marked
point $p$ and arbitrary local orientation at $p$;
\\(iv) $(\iota_A(\alpha'),\alpha'')_A=\Phi_{S^2,p',p''(-)}(\alpha',\alpha'')$ where
$S^2$ is a sphere with marked points $p'$ and $p''$, local
orientation at $p'$ coincides with a fixed global orientation of
$S^2$  and local orientation at $p''$ does not coincide with the
fixed global orientation of $S^2$;
\\(v)  $(\iota_B(\beta'),\beta'')_B=\Phi_{D^2,q',q''(-)}(\beta',\beta'')$ where
$D^2$ is a disc with marked   points $q'$ and $q''$, local
orientation at $q'$ coincides with a fixed global orientation of
$D^2$  and local orientation at $q''$ does not coincide with the
fixed global orientation of $D^2$.
\\(vi) $(K_A^\otimes,\alpha'\otimes\alpha'')_A=(\alpha',\alpha'')_A$  for any
$\alpha',\alpha''\in A$;
\\(vii) $(K_B^\otimes,\beta'\otimes\beta'')_B=(\beta',\beta'')_B$  for any
$\beta',\beta''\in B$;
\end{lemma}

\proof The proof is straightforward $\Box$

\subsection{Klein Topological Field Theory of $G$-coverings}
\label{klein_g} Let $G$ be a finite group and $K$ be a subgroup of
$G$ such that $\cap_{g\in G}gKg^{-1}=\{e\}$. We shell construct a
KTFT corresponding to $G$-coverings of surfaces by seamed surfaces
such that stabilizers of all simple points of cover seamed surface
are conjugated to $K$.

Denote by $\mathcal{A}$ the set of all possible values of local
invariants of $G$-coverings $T\to\Omega$ at interior marked points
of the base. The set $\mathcal{A}$ can be identified with the set of
conjugacy classes of elements of the factor-group $N=N_G(K)/K$
(see section \ref{local}). Denote by $A$ a vector space of formal
linear combinations of elements of $\mathcal{A}$.

Denote by $\mathcal{B}$ the set of all possible values of local
invariants of $G$-coverings $T\to\Omega$ at boundary marked points
of the base. The set $\mathcal{B}$ can be identified with the set of
equivalence classes of ordered pairs of subgroups $(S_1,S_2)$ of the
group $G$ such that $K\subset S_1\cap S_2$; the equivalence is
defined by formula  $(S_1,S_2)\sim (gS_1g^{-1},gS_2g^{-1})$ for any
$g\in N_G(K)$ (see section \ref{local}). Denote by $B$ a vector
space of formal linear combinations of elements of $\mathcal{B}$.

For any surface $\Omega$ Hurwitz numbers
$\mathcal{H}(\Omega;\mathfrak{a},\mathfrak{b})$ where
$\mathfrak{a}:\Omega_a\to\mathcal{A}$ and
$\mathfrak{b}:\Omega_b\to\mathcal{B}$ are maps used in the
definition of Hurwitz numbers, induce the linear function
$\mathcal{H}_\Omega:V_\Omega\to \mathbb{C}$.

We will equip vector spaces $A$ and $B$ by a set of data
needed to define a Klein functor (section \ref{klein_f}).

By $1_A$ denote the conjugacy class of the unit in group
$N=N_G(K)/K$; by definition, $1_A\in A$.

Let $\mathcal{R}$ be the set of equivalent classes of
subgroups $S\subset G$ such that $S\supset K$; the equivalence is
defined by formula $S\sim nSn^{-1}$ for any $n\in N_G(K)$.
By $1_B$ denote the sum $\sum_{[S]\in\mathcal{S}}[(S,S)]$ where $S$
is a representative of equivalence class $[S]\in\mathcal{R}$;
by definition, $1_B\in B$.

Denote by  $\iota_A:A\to A$ and $\iota_B:B\to B$  involute linear
maps induced by involute maps $\alpha\to\alpha^*$ and
$\beta\to\beta^*$ of $\mathcal{A}$ and $\mathcal{B}$ respectively
(see section \ref{local}).

By $K_A^\otimes$ and $K_B^\otimes$ denote the elements
$K_A^\otimes=\sum_{\alpha\in\mathcal{A}}|\Aut\alpha|\alpha\otimes\alpha^*$
and $K_B^\otimes=\sum_{\beta\in\mathcal{B}}|\Aut\beta|\beta\otimes\beta^*$
respectively; by definition, $K_A^\otimes\in A\otimes A$,
$K_B^\otimes \in B\otimes B$.

Denote by  $U$ the element  $U=\sum_{g\in N}g^2$ of $A$.

Thus, we defined the set of data
\\$\mathcal{N}=\{A, B, 1_A\in A, 1_B\in B, K_A^\otimes\in A\otimes A, K_B^\otimes\in B\otimes B, U\in
A,\iota_A,\iota_B\}$, which is needed to construct a Klein functor  (see section
\ref{klein_f}).

\begin{theorem} The system of linear functionals
$\mathcal{H}=\{\mathcal{H}_\Omega\}$ is a Klein topological field
theory. We call it \textit{topological field theory of
$G$-coverings of surfaces by seamed surfaces}.
\end{theorem}

\proof

Let $T\to\Omega$ be a $G$-covering of $\Omega$ and $\alpha\in\Omega$
be a local invariant at interior marked point $p\in \Omega$. Note
that local orientation at $p$ is used to define $\alpha$ (see
section \ref{local}). Changing the local orientation at $p$ leads to
replacing a conjugacy class  $\alpha\subset N_G(K)/K$ by the conjugacy class
by $\alpha^*=[a^{-1}]$ where $a\in\alpha$.
That relation is equivalent to axiom $2^\circ$ of KTFT for interior marked points.
Axiom $2^\circ$ for boundary marked points can be verified similarly.

Axiom $1^\circ$ is a corollary of the topological invariance of Hurwitz numbers with
respect to isomorphisms of surfaces preserving local orientations at marked points.

Axiom $3^\circ$ (nondegeneracy of the form) follows from the straightforward
calculations of Hurwitz numbers for a sphere with two interior
marked points and a disc with two boundary marked points.

Axiom $4^\circ$ is evidently satisfied for the marking off a new
interior point $p$ with $\mathfrak{a}(p)=1_A$ because the conditions
of general position are actually required from a  covering at that
point. The same is true for boundary marked points too.

Axiom $5^\circ$ follows from lemma \ref{lemma_cut}.

Axiom $6^\circ$ follows from the definition of Hurwitz numbers.

$\Box$

\section{Cardy-Frobenius algebras associated with actions of groups}
\label{cardy}

It follows from \cite{AN} that KTFTs up to the equivalence are in
one-to-one correspondence with a certain class of algebras
up to isomorphisms. Algebras of that class was called
'structure algebras' in \cite{AN} and renamed to 'Cardy-Frobenius
algebras' in \cite{AN2}. We use the latter term here.

\subsection{Definition of Cardy-Frobenius algebra}
\label{definition_c}

A Frobenius algebra $D$ is an associative algebra such that there exists a
nondegenerate symmetric invariant bilinear form $(x',x'')_D$ on $D$.
The invariance means $(x'x'',x''')_D=(x',x''x''')_D$.

A finite dimensional Frobenius algebra $D$ with unit $1_D\in D$
endowed with an involute anti-automorphism $d\mapsto d^*$ (i.e.
$(x'x'')^*=x''^*x'^*$) and a linear form $l_D(x)$ such that the
bilinear form $(x',x'')_D=l_D(x'x'')$ is symmetric, nondegenerate
and $l_D(x^*)=l_D(x)$ is called {\it an equipped Frobenius algebra}.

Let $D$ be an equipped  Frobenius algebra. Choose a basis
$\alpha_1,\dots,\alpha_n$ of vector space $D$. Denote by
$((F_{\alpha_i,\alpha_j}))$ the matrix of the bilinear form,
$F_{\alpha_i,\alpha_j}= (\alpha_i,\alpha_j)_D$, by
$((F^{\alpha_i,\alpha_j}))$ the matrix of the dual form and by
$(({F^*}_{\alpha_i,\alpha_j}))$ (resp.,
$(({F^*}^{\alpha_i,\alpha_j}))$) the matrix of twisted bilinear form
${F^*}_{\alpha_i,\alpha_j}=(\alpha_i,\alpha_j^*)_D$ (resp. dual to
twisted bilinear form).

An element $K_{D}=F^{\alpha_i\alpha_j}\alpha_i\alpha_j$ of $D$ is called
\textit{a Casimir element}. In that formula (and
throughout the paper) the sum over repeated indexes is assumed.
Define also  \textit{twisted Casimir element} $K_{D}^*\in D$ by
formula $K_{D}^*= {F^*}^{\alpha_i\alpha_j}\alpha_i\alpha_j$. It can
be shown that Casimir element and twisted Casimir element don't
depend on the choice of basis.

Fix a set of data $\mathcal{H}= \{A,B,\phi, U\}$, where $A$ and $B$
are equipped Frobenius algebras, $A$ is a commutative algebra,
$\phi:A\to B$ is a homomorphism of algebra $A$ into the center of
algebra $B$, $U\in A$ is an element. Denote by $H$ an algebra
structure on $A\oplus B$ defined by multiplications in $A$ and $B$
and by formula $ab=ba=\phi(a)b$ for $a\in A, b\in B$. Thus, $A$ is a
subalgebra and $B$ is two-sided ideal of $H$.

The algebra $H=A\oplus B$ together with the set of data $\mathcal{H}=
\{A,B,\phi, U\}$ is called {\it a
Cardy-Frobenius algebra} if the following conditions hold:

(1) $\phi(x^*)=(\phi(x))^*$

(2) (Cardy condition) For any $x,y\in B$ denote by $W_{x,y}:B\to B$
a linear operator $W_{x,y}(z)=xzy$. Let $\phi^*: B\rightarrow A$ be
the linear map dual to $\phi: A\to B$ (i.e. $(a,\phi^*(b))_A=
(\phi(a),b)_B$, for $a\in A, b\in B$). It is required that for any
$x,y\in B$ the following identity holds:

$$(\phi^*(x),\phi^*(y))_A= trW_{x, y}$$

(3) $U^2=K_A^*$ and $\phi(U)=K_B^*$

\subsection{Semisimple Cardy-Frobenius algebras}
\label{semisimple}

Direct sums of Cardy-Frobenius algebras, ideals etc. can be defined in usual way.
We skip the details.

\begin{example}\label{ex1}
Let $A$ be an algebra of complex numbers $\mathbb{C}$
equipped with identical involution and linear form $l_A(z)=\mu^2 z$
where $\mu\in\mathbb{C}$ is nonzero number. Then $A$ is equipped
Frobenius algebra.

By definition, put $U=\frac{1}{\mu}\in A$.

Let $B$ be a matrix algebra $\mathbb{M}(n,\mathbb{C})$ equipped with
involute anti-automorphism $X\mapsto X'$, where $X'$ is transposed
matrix, and linear form $l_B(X)=\mu \tr X$. Then  $B$ is equipped
Frobenius algebra. Define a homomorphism $\phi:A\to B$ by the
equality $\phi(1)=E$ where $E\in B$ is the identity matrix.

It can be shown that the set of data $\{A,B,\phi,U\}$ defines the simple
Cardy-Frobenius algebra $H=A\oplus B$.
We call $H$ a Cardy-Frobenius algebra $H_{(n,\mu)}^+$

\end{example}

\begin{example} Let $A$ be equipped Frobenius algebra from previous
example. By definition, put $U=-\frac{1}{\mu}\in A$.

Let $B$ be a matrix algebra $\mathbb{M}(2m,\mathbb{C})$ of even
order. A matrix $X\in B$ we may present in block form as $X=
\left(%
\begin{array}{cc}
  m_{11}&  m_{12}
   \\
  m_{21}& m_{22} \\
\end{array}%
\right)$. Define the involute anti-automorphism $X\mapsto X^\tau$ by
formula
 $X^\tau=
\left(%
\begin{array}{cc}
  m_{22}'& -m_{12}'
   \\
 -m_{21}'& m_{11}' \\
\end{array}%
\right),$ where $m'$ denotes transposed matrix. Then algebra $B$
equipped with involute anti-automorphism $X\mapsto X^\tau$ and
linear form $l_B(X)=\mu \tr X$ is equipped Frobenius algebra.

Define homomorphism $\phi:A\to B$ by the equality $\phi(1)=E$ where $E\in B$ is the
identity matrix.

It can be shown that the set of data $\{A,B,\phi,U\}$ defines the simple
Cardy-Frobenius algebra $H=A\oplus B$. We call $H$ a Cardy-Frobenius
algebra $H_{(2m,\mu)}^-$
\end{example}

\begin{example}\label{ex3}
Let $A=\mathbb{C}\oplus\mathbb{C}$. Define an involution by
the equality $(x,y)^*=(y,x)$ for $(x,y)\in \mathbb{C}\oplus\mathbb{C}$
and the linear form by formula $l_A(x,y)=\mu^2 (x+y)$. Clearly, $A$
is commutative equipped Frobenius algebra. Put $U=0$.

Let $B=\mathbb{M}(n,\mathbb{C})\oplus\mathbb{M}(n,\mathbb{C})$ be
the direct sum of two matrix algebras of equal order. Then algebra
$B$ equipped with the involute anti-automorphism $(X,Y)\mapsto
(Y',X')$ and linear form $$l_B(X,Y)=\mu (\tr X + \tr Y)$$ is an
equipped Frobenius algebra.

Define a homomorphism $\phi:A\to B$ by the equality $\phi(x,y)=(xE,yE)$.

It can be shown that the set of data $\{A,B,\phi,U\}$ defines the simple
Cardy-Frobenius algebra $H=A\oplus B$. We call $H$ a Cardy-Frobenius
algebra $H_{(2n,\mu)}^0$

\end{example}

\begin{example}

Denote by $A_{\mu}$ the pair $(A,U)$ from example \ref{ex1}. Denote
by $A_{\mu}^0$ the pair $(A,U)$ from example \ref{ex3}. Note that
algebras $A_{\mu}\oplus\{0\}$ and $A_{\mu}^0\oplus\{0\}$ satisfy all
conditions of Cardy-Frobenius algebras except one: $B$ has no unit.

\end{example}

It was proven \cite{AN}, that any Cardy-Frobenius algebra such that
both algebras $A$ and $B$ are semisimple is isomorphic to the direct
sum of simple Cardy-Frobenius algebras of types $H_{(n,\mu)}^+$,
$H_{(2m,\mu)}^-$, $H_{(2n,\mu)}^0$ and commutative equipped
Frobenius algebras with fixed $U\in A$ of types $A_\mu$ and
$A_\mu^0$.

\subsection{Cardy-Frobenius algebras associated with actions of groups}
\label{equipped} We will construct  here a wide class of
Cardy-Frobenius algebras $H=A\oplus B$. Among algebras $B$
(noncommutative parts  of Cardy-Frobenius algebras) there are, for
example, all Hecke algebras of finite groups \cite{S}. In the next
section we will show that one of Cardy-Frobenius algebras of that
class coincides with Cardy-Frobenius algebra of $G$-coverings of
surfaces by seamed surfaces.

Let $N$ be a finite group. Group algebra $\mathbb{C}[N]$ can be
equipped with involute antiautomorphism $x\mapsto x^*$ defined on
$g\in N$  by formula $g^*=g^{-1}$ and linear form $l(x)$ defined by
formula $l(g)=\delta_{g,e}$ for $g\in N$ ($e$  denotes the unit of
$N$).

Denote by $A$ the center of the group algebra . Clearly, $A$ is invariant under the involution
$x\mapsto x^*$. It can be easily shown that the algebra $A$ endowed with that involution
and the restriction of linear form $l(x)$ is commutative equipped Frobenius algebra.

Denote by $U$ the element $U=\sum_{g\in N} g^2$ of group algebra. Clearly, $U\in A$.

\begin{lemma} The square of element $U$ is equal to the twisted Cazimir
element: $U^2=K_A^*$
\end{lemma}
\proof We have  $U^2=\sum_{g,h\in N} g^2h^2$.

Denote by $\mathcal{A}$ the set of conjugacy classes of the group $N$.
For a conjugacy class $\alpha$ put
$E_\alpha=\sum_{g\in \alpha}g\in\mathbb{C}[N]$. It is known that  elements $E_\alpha$
belong to the center $A$ of the group algebra and form a basis of $A$.

We have  $K_A^*=\sum_{\alpha\in\mathcal{A}}|\Aut\alpha|E_\alpha^2=
\sum_{\alpha\in\mathcal{A}}|\Aut\alpha|(\sum_{g\in\alpha}g)^2=
\sum_{g,h\in N} ghg^{-1}h$. Using substitutions $a=gh$ and $b=h^{-1}g^{-1}h$ we get
$\sum_{g,h\in N} ghg^{-1}h=\sum_{a,b\in N} a^2b^2$. Therefore, $U^2=K_A^*$ $\Box$

Suppose, finite group $N$ acts effectively on a finite set $X$. Define by $X^n$ the
direct product of $n$ copies of $X$. Group $N$ acts on $X^n$ by
formula $g(x_1,\dots,x_n)= (g(x_1),\dots,g(x_n))$ where $x_1,\dots,x_n\in X$, $g\in N$.
Denote by $\mathcal{B}_n$ the set of orbits $X^n/N$.

Denote by $\Aut\bar x$ the stabilizer of element $\bar x=(x_1,\dots,x_n)\in X^n$
Clearly, $\Aut\bar x=\cap_i\Aut x_i$.

Let $B$ be the vector space of formal linear combinations with complex
coefficients of elements of $\mathcal{B}_2$.

Fix an ordered set $(\beta_1,\dots,\beta_n)$ of elements $\beta_i\in\mathcal{B}_2$.
By definition, an element $\beta_i$ is an orbit of the action of group $N$ on $X^2$.
Denote by $\mathcal{B}_n(\beta_1,\dots,\beta_n)$ the set of orbits $[\xi]$ in $X^n$
such that the following condition holds: for a representative
$\xi=(x_1,\dots,x_n)\in [\xi]$ the pair
$(x_1,x_2)$ belongs to the orbit $\beta_1$,
$(x_2,x_3)$ belongs to $\beta_2$, $\dots$, $(x_{n-1},x_n)$ belongs to $\beta_{n-1}$,
$(x_n,x_1)$ belongs to $\beta_n$.

Define a number $T_{\beta_1,\dots,\beta_n}$ by formula
$$T_{\beta_1,\dots,\beta_n}=\sum_{[\xi]\in\mathcal{B}_n
(\beta_1,\dots,\beta_n)}\frac{1}{\Aut\xi},$$ where $\xi$ is a
representative of an orbit $[\xi]\in\mathcal{B}_n(\beta_1,\dots,\beta_n)$.

Numbers  $\{T_{\beta_1,\dots,\beta_n}\}$ can be considered as a tensor, it corresponds
to a polylinear form $B^{\otimes n}\to \mathbb{C}$.

The correspondence $(x_1,x_2)\mapsto(x_2,x_1)$ generates an involute
linear transformation of $B$ which we denote by
$\beta\mapsto\beta^*$.

Tensor $T_{\beta_1,\beta_2}$ defines a bilinear form on $B$, which we denote
by $(\beta',\beta'')_B$. By direct calculations we obtain
$(\beta_1,\beta_2)_B=\delta_{\beta_1,\beta_2^*}\frac{1}{|\Aut\beta_1|}$. Therefore,
bilinear form $(\beta',\beta'')_B$ is symmetric and non-degenerate.

Using the bilinear form $(\beta',\beta'')_B$ we may convert tensor $T_{\beta_1,\beta_2,\beta_3}$
to a multiplication on $B$. Indeed, formula $(\beta_1\cdot\beta_2,\beta_3)=T_{\beta_1,\beta_2,\beta_3}$
correctly defines an element $\beta_1\cdot\beta_2$. The obtained algebra is denoted by the same letter $B$.

Denote by $\mathcal{B}_2^circ$ the set of orbits of elements $(x,x)\in X^2$, $x\in X$.
By straightforward computations we get that the element
$1_B=\sum_{\beta\in\mathcal{B}_2^\circ}\beta$ is a unit of algebra $B$.
Formula $l_B(\beta)=(\beta,1_B)_B$ defines a linear form.

\begin{lemma} \label{lemma_A}
Algebra $B$ endowed  with linear form  $l_B$ and involution
$\beta\mapsto\beta^*$ is an equipped Frobenius algebra.
\end{lemma}

The proof is by direct calculations $\Box$

\smallskip

Denote by $V_X$ the  vector space of formal linear
combinations with complex coefficients of elements of set $X$ and extend the action of $N$ on
$X$ to the representation $\rho$  of $N$ in $V_X$.

For any $\beta\in\mathcal{B}_2$ define the operator $V_\beta\in \End V_X$ by
formula $V_\beta=\sum_{(x_1,x_2)\in\beta}E_{x_1,x_2}$. Here
$E_{x_1,x_2}$ is the elementary matrix:
$E_{x_1,x_2}x=\delta_{x_2,x}x_1$ for any $x\in X$.
Denote by $\nu$ the linear map $B\to\End V_X$,
induced by the correspondence $\beta\to V_\beta$.

\begin{lemma} \label{centralizer}
The map $\nu:B\to\End V_X$ is a faithful representation of algebra $B$.
The image $\nu(B)\subset \End V_X$ coincides with the centralizer of
the image $\rho(N)\subset\End V_X$ of group $N$ under the representation $\rho:N\to\End V_X$.
\end{lemma}

\proof  Let $((M_{x,y}))$ be the matrix of a linear transformation $M:V_X\to V_X$
in the basis, consisting of elements $x\in X$. Then by direct calculations we obtain that
$F$ commutes with $\rho(N)$ if and only if matrix elements  $M_{x,y}$ and $M_{g(x),g(y)}$ are equal for any
$(x,y)\in X^2$ and $g\in N$. This claim is equivalent to the claim of lemma $\Box$

\begin{corollary}
Algebra $B$ is semisimple equipped Frobenius algebra.
\end{corollary}

\proof Indeed, the centralizer of a semisimple subalgebra of semisimple algebra is
semisimple $\Box$

The constructed equipped Frobenius algebra $B$ depends on group $N$
and its action on finite set $X$. We denote that algebra by
$B_{N:X}$. Algebra $A$ acts on $V_X$ because $N$ acts on $V_X$.
Clearly, that representation of $A$ is faithful. By lemma
\ref{centralizer}, the image $\rho(A)$ is contained in the center of
the image of algebra $B$. The representation of $B$ in $V_X$ is
faithful, therefore, we obtain the homomorphism  of $A$ into $B$.
That homomorphism we denote by $\phi$.

\begin{theorem} Let a finite group $N$ acts effectively on a set $X$. Then
defined above objects: the algebras $A$, $B$, the homomorphism
$\phi:A\to B$ and the element $U\in A$, generate a Cardy-Frobenius
algebra $H=A\oplus B$.
\end{theorem}

\proof Clearly, that involute anti-automorphism
$\beta\mapsto\beta^*$ acts on the image of $B$ in $\End V_X$ as the
transposition $M\mapsto M'$ of matrices. Denote by $V_\alpha $ the
image of element $E_\alpha\in A$ in $\End V_X$ under the
representation $\rho:\mathbb{C}[N]\to\End V_X$. Obviously,
$V_\alpha'=V_{\alpha^*}$. Therefore, condition (2) of the definition
of Cardy-Frobenius algebra is satisfied.

The proof of Cardy relation is literally the same as the
proof of that relation in \cite{AN} for special class of
Cardy-Frobenius algebras, namely, Cardy-Frobenius algebras
associated with finite groups. The proof is rather
technical, we outline here main steps.

\bigskip
Compute first right side expression $R=\tr W_{x,y}$.

(1) Scalar product of elements $\beta_i\in\mathcal{B}_2$  can be computed in the representation
of algebra $B$ in $\End {V_X}$ as follows:
$(\beta_1,\beta_2)_B=\frac{1}{|G|}\tr ( V_{\beta_1}V_{\beta_2})$ where $V_{\beta_i}$ are operators
defined above.

(2) By direct calculations using exact formulas for matrix elements of $V_\beta$ we get
$R=\sum_\beta \frac{|\Aut \beta|}{|G|}\tr (V_\beta V_{\beta_1}V_{\beta^*}V_{\beta_2})$.

(3) Multiplying matrices in (2) we obtain

$R=\frac{1}{|G|} |P|$, where $P$ is the set of tuples
$(x_1,x_2,x_3,x_4, g)$, $x_i\in X$, $g\in G$ such that
$(x_1,x_2)\in\beta_1$, $(x_3,x_4)\in\beta_2$ and $g(x_2)=x_1$,
$g(x_3)=x_4$.

\bigskip

Compute left side expression $L=(\phi^*(\beta_1),\phi^*(\beta_2))_A$.

For any $\alpha\in\mathcal{A}$ we have $V_\alpha=\sum_{x\in X, g\in
\alpha} E_{g(x),x}$, where $E_{x_1,x_2}$ is an elementary matrix.
(Recall that $V_\alpha \in\End V_X$ is the image of element
$E_\alpha\in A$ in $\End V_X$ under the representation
$\rho:\mathbb{C}[G]\to\End V_X$).

(4) By direct calculations we obtain
$L=\sum_\alpha |\Aut\alpha| \frac{1}{|G|^2}\tr (V_\alpha V_{\beta_1}) \tr (V_{\alpha^*} V_{\beta_2})$.

(5) The equality (4) can be rewritten in the form

$L=\sum_\alpha |\Aut\alpha| \frac{1}{|G|^2} |Q|$
where $Q$ is the set of tuples
$(x_1,x_2,x_3,x_4, a_1, a_2)$, $x_i\in X$, $a_j\in \alpha$ such that
$(x_1,x_2)\in\beta_1$, $(x_3,x_4)\in\beta_2$ and  $a_1(x_1)=x_2$, $a_2(x_4)=x_3$.

(6) Representing element $a_2$ in the form $ga_1g^{-1}$ we can rewrite (5) as follows:

$L=\sum_\alpha \frac{1}{|G|^2}|Q'|$ where $Q'$ is the set of tuples
$(x_1,x_2,x_3,x_4, a, g)$, $x_i\in X$, $a\in \alpha$, $g\in G$ such that
$(x_1,x_2)\in\beta_1$, $(x_3,x_4)\in\beta_2$ and $a(x_1)=x_2$, $gag^{-1}(x_4)=x_3$.

(7) The equality $gag^{-1}(x_4)=x_3$ can be rewritten as $a(x_4')=x_3'$ where $x_i'=g^{-1}(x_i)$.
For two fixed $g=g_1$ and $g=g_2$ the number of tuples $(x_1,x_2,x_3,x_4, a, g_1)\in Q'$
is equal to the number of tuples $(x_1,x_2,x_3,x_4, a, g_2)\in Q'$
since the pairs $(x_3,x_4)$ run throughout all representatives of class $\beta_2$.
Therefore, $L=\frac{1}{|G|} |P|$. Thus, we prove $L=R$.

\bigskip

To complete the proof, it is necessary  to check that $\rho(U)=\nu(K_B^*)$.
By definition, $U=\sum_{g\in N}g^2$. Therefore,
$\rho(U)= \sum_{g\in N}E_{g^2(x),x}$. Denote $\nu(K_B^*)$ by $K$.
By direct calculations we  we obtain
$K=\sum_\beta|\Aut\beta|\sum_{(x_1,x_2)\in\beta}E_{x_1,x_2}\sum_{(y_1,y_2)\in\beta}E_{y_1,y_2}$.
We can substitute $(y_1,y_2)$ by $(g(x_1),g(x_2)$ for some element $g\in G$  because both pairs
$(x_1,x_2)$ and $(y_1,y_2)$ belong to the same class $\beta$. Taking in account the substitution
and multiplying matrices  we get $K=\sum_{x\in X, g\in G} E_{x,g^2(x)}$.
Thus, $K=\rho(U)$ $\Box$

\bigskip

We denote constructed Cardy-Frobenius algebra by $H_{G:X}$.

\begin{example} Let symmetric group $S_n$
acts on a finite set $X$ of $n$ elements. Then algebra $B_{S_n:X}$
is isomorphic to the algebra of bipartite graphs \cite{AN1,AN2}.
Moreover, $H_{S_n:X}$ coincides with  Cardy-Frobenius algebra of
$n$-coverings of Klein surfaces by seamed surfaces \cite{AN1,AN2}.
\end{example}

\begin{example} Let $X$ be the set of involutions in a finite group $G$.
Group $G$ acts on $X$ by conjugations. Then $H_{G:X}$
coincides with constructed in \cite{AN} Cardy-Frobenius algebra associated with
$G$. For symmetric group $G=S_n$ the latter algebra coincides with Cardy-Frobenius
algebra of coverings of degree $n$ for real algebraic curves \cite{AN}.
\end{example}

\begin{example} Let $S$ be a subgroup of finite group $G$ and $X=G/S$.
Then algebra $B_{G:X}$ is isomorphic to Hecke algebra. Indeed, pairs of cosets
$(xB,yB)\in X\times X$ are in one-to-one correspondence with
double-cosets $BgB$ and formulas for the multiplication in algebra
$B_{G:X}$ coincide with formulas for the multiplication of double-cosets
in Hecke algebra \cite{S}.
\end{example}

\section{Formula for Hurwitz numbers of $G$-coverings of surfaces
by seamed surfaces} \label{formula}

\subsection{Cardy-Frobenius algebra corresponding to a Klein topological field theory}
\label{cardy_a}

Klein topological field theories (KTFTs) are in one-to-one correspondence with  Cardy-Frobenius
algebras \cite{AN}. Following \cite{AN} we briefly outline the construction of Cardy-Frobenius algebra
corresponding to a KTFT.

By definition, Cardy-Frobenius algebra is generated by the following set
of objects: a commutative equipped Frobenius algebra $A$,
an equipped Frobenius algebra $B$, a homomorphism $\phi:A\to B$,
an element $U\in A$.

Let $\mathcal {F}$ be a KTFT. By definition, $\mathcal {F}$ is a
system of linear forms $\{\Phi_{\Omega}:V_\Omega\to\mathbb{C}\}$,
defined for each surface $\Omega$ and satisfying certain axioms (see
section \ref{klein_t}). Vector spaces $V_\Omega$ are defined by the
use of fixed Klein functor $\mathcal{V}$. In turn, Klein functor
corresponds to a set of data $\mathcal{N}$; $\mathcal{N}$ includes,
among others, two vector spaces $A$ and $B$, elements $1_A\in A$,
$U\in A$, $1_B\in B$, involute linear maps $\iota_A:A\to A$,
$\iota_B:B\to B$ (see section \ref{klein_f}). Using Klein functor
$\mathcal{V}$ and KTFT $\mathcal{F}$, one can define  algebra
structures on vector spaces $A$ and $B$ and equip those algebras
with required in the definition of Cardy-Frobenius algebras
additional objects.

Let $\Omega$ be a sphere with two marked points,
$\Omega=(S^2,p',p'')$. Then corresponding to $\Omega$ linear form
$\Phi_{(S^2,p',p'')}:A\otimes A\to\mathbb{C}$  induces bilinear form
$(x,y)_A$ on $A$. By axioms, the form is symmetric and
nondegenerate. Similarly, disc $(D^2,q',q'')$ with to boundary
marked points give rise to nondegenerate bilinear form on $B$.

By definition, tri-linear map $\Phi:A\otimes A\otimes
A\to\mathbb{C}$ corresponds to a sphere with three marked points.
Using bilinear form on $A$,  a tri-linear form can be converted into
the multiplication $A\otimes A\to A$. Axioms of KTFT provide the
commutativity and associativity of that multiplication and the
invariance of the bilinear form. By construction, algebra $A$ is
equipped with the involute map $\iota_A$, fixed element $1_A$ and
linear form $l_A(x)=(x,1_A)_A$. It can be shown that $A$ is
commutative equipped Frobenius algebra.

Similarly, a triangle, i.e. a disk with three marked boundary
points, gives rise to algebra structure on $B$. It can be
proven that $B$ is (typically, noncommutative) equipped Frobenius algebra.

Let $\Omega=(D^2,p,q)$ be a disc with one boundary and one interior marked points.  Then
linear form $\Phi_{(D^2,p,q)}:A\otimes B\to\mathbb{C}$ corresponds to $\Omega$.
Using bilinear forms, the map $\Phi_{(D^2,p,q)}$ can be transformed into the
map  $\phi:A\to B$. It can be shown that $\phi$ is a homomorphism of algebras.

In \cite{AN} was proven that the data $\mathcal{H}=\{A,B,\phi,U\}$
correctly defines the structure of  Cardy-Frobenius algebra
on $H=A\oplus B$. In the proof different systems of
simple cuts for certain surfaces were considered.
For example, two essentially different systems of
simple cuts of a cylinder with two marked boundary points belonging to different
connected components of the boundary were used to prove Cardy relation.

\subsection{Cardy-Frobenius algebra of $G$-coverings of surfaces by seamed surfaces}
\label{cardy_g}

Let $G$ be a finite group, $K\subset G$ be a subgroup such that
$\cap_{g\in G}gKg^{-1}=\{e\}$. In section \ref{klein_g} KTFT
$\mathcal{H}$ of $G$-coverings with stabilizers of simple points
conjugated to $K$  of surfaces by seamed surfaces was associated
with the pair $G\supset K$. Denote by $H$ the Cardy-Frobenius
algebra corresponding to KTFT $\mathcal{H}$.

Denote by $N$ the factor group $N_G(K)/K$ and by
$X$ the set of all subgroups $S\subset G$ such that $S\supset K$.
Conjugations $S\to hSh^{-1}$, $h\in N_G(K)$ generate the action of
$N$ on $X$ because elements of $K$ acts trivially.
Denote by $H_{N:X}$  Cardy-Frobenius algebra associated with that action
(see section \ref{equipped}).

\begin{theorem} Cardy-Frobenius algebra $H_{N:X}$ is isomorphic
to  Cardy-Frobenius algebra $H$ of $G$-coverings with stabilizers of simple
points conjugated to $K$ of surfaces by seamed surfaces.
\end{theorem}

\proof Both algebras $H$ and $H_{N:X}$ were defined by choosing
bases of vector spaces and explicit formulas for the multiplications
of elements of the bases.  Comparing those data for $H$
(see section \ref{cardy_a}) and $H_{N:X}$ (see section \ref{equipped})
we observe that bases of both Cardy-Frobenius algebras are
in one-to-one correspondence and formulas for the multiplications coincide.
Thus, the algebras are isomorphic. We skip the details $\Box$

\begin{example} Let $G$ be the smallest finite
simple group $A_5$ and $K$ be a subgroup of two elements,
$K=Z_2$. Denote by $H=A\oplus B$ the Cardy-Frobenius algebra of
the KTFT of $A_5$-coverings with stabilizers of simple points conjugated
to $Z_2$ of surfaces by seamed surfaces. Clearly, subgroup $N_G(K)$
is isomorphic to $Z_2\times Z_2$, therefore, the group $N=N_G(K)/K$
consists of two elements.

We shell computed $H$
using the isomorphism $H\approx H_{N:X}$ where $X$ is the set of
all subgroups of $A_5$ containing $K$.

By the above, algebra $A$ is isomorphic to the center of
group algebra $\mathbb{C}[N]$. Therefore,
$A=\mathbb{C}\oplus\mathbb{C}$.

We claim that $B=M(1,\mathbb{C})\oplus M(5,\mathbb{C})$.
Indeed, there is six subgroups containing $K$: (1) $K=Z_2$,
(2) $D_4=Z_2\times Z_2$, (3) $A_4$,
(4) $G=A_5$ and (5-6) two conjugated in $G$ subgroups of order $10$; those
subgroups are denoted by $L_{10}'$ and $L_{10}''$, each of them is semidirect
product of $Z_2$ and $Z_5$ with non-trivial action of $Z_2$ on $Z_5$.
Thus, the set $X$ consists of six elements.
Non-unit element $n$ of group $N$ acts of $X$ by permutation of  $L_{10}'$
and $L_{10}''$, other elements of $X$ are fixed points of $n$.

Denote by $V_X$ the vector space of linear combinations with complex coefficients
of elements of $X$. Clearly, the centralizer of the image of $N$ in $\End V_X$ is
isomorphic to $M(1,\mathbb{C})\oplus M(5,\mathbb{C})$. By lemma
\ref{centralizer} $B$ is isomorphic to that centralizer.
Homomorphism $\phi:A \to B$ is actually an isomorphism of $A$ with
the center of algebra $B$.

Equipments of algebras $A$ and $B$ can be easily computed by lemma \ref{lemma_A}.
Thus, Cardy-Frobenius algebra $H$ is described.

\end{example}

Note, that $G$-coverings with certain restrictions on local topological
invariants at marked points can be described by a similar approach.
Indeed, let $G$ be a finite group, $K$ be its subgroup and
$\Sigma$ be a set of subgroups such that for any $S\in \Sigma$
(1) $S\supset K$ and (2)  $gSg^{-1}\in\Sigma$ for any $g\in N_G(K)$.
Then group $N=N_G(K)/K$ acts on $\Sigma$. Denote by
$H_{N:\Sigma}$ a Cardy-Frobenius algebra corresponding to that action.
Clearly, $H_{N:\Sigma}$ is a subalgebra of Cardy-Frobenius
algebra $H_{N:X}$ where $X$ is the set of all subgroups containing $K$.
The subalgebra $H_{N:\Sigma}$can be considered as
Cardy-Frobenius algebra of $G$-coverings of surfaces by seamed surfaces
with restricted sets of local invariants.

\begin{example}
Let $G$ be symmetric group $S_n$, $K=\{e\}$ and $\Sigma$ the set of subgroups of order $2$
and $1$. Then $N=N_G(\{e\})$ coincides with $G$. It can be shown that
Cardy-Frobenius algebra $H_{S_n:\Sigma}$ is isomorphic to
Cardy-Frobenius algebra of real algebraic curves constructed in \cite{AN}.
\end{example}

\subsection{Formula for Hurwitz numbers} \label{formulahurw}

Let $G$ be a finite group, $K\subset G$ be a subgroup and $H$ be
Cardy-Frobenius algebra corresponding to KTFT of $G$-coverings of surfaces by seamed surfaces
such that all stabilizers of simple points of cover space are conjugated to $K$.
Following \cite{AN}, we will express Hurwitz numbers
$\mathcal{H}(\Omega,\alpha,\beta)$ in terms of algebra $H$.

First, we may restrict ourselves to connected surfaces $\Omega$ due to multiplicativity
of Hurwitz numbers with respect to the disjoint union of surfaces.

Second, a connected surface is defined, up to isomorphism, by the following topological
invariants:
\begin{itemize}
\item orientability $\varepsilon$ ($\varepsilon=1$ for orientable surfaces, $\varepsilon=0$ for
non-orientable surfaces;
\item genus $g$; for non-orientable surfaces $g$ is half-integer: $g=\frac{1}{2}$ for projective
plan, $g=1$ for Klein bottle etc.
\item number $s$ of connected components of the boundary $\partial\Omega$;
\item number $m$ of interior marked points;
\item numbers $m_1,\dots,m_s$ of marked points at each boundary contour.
\end{itemize}

If $\Omega$ is orientable we choose local orientations at all marked points compatible with a
fixed  orientation of $\Omega$. If $\Omega$ is non-orientable then we choose arbitrary local orientations
at interior boundary points and for each boundary contour we choose local orientations at
boundary marked points lying at that contour compatible with fixed orientation of the contour.

Third, for given map $\alpha:\Omega_a\to\mathcal{A}$ denote by $\alpha_1,\dots\alpha_m$ images
of $\alpha$, i.e. interior fields. For given map $\beta:\Omega_b\to\mathcal{B}$
denote by $\beta^i_1,\dots,\beta^i_{m_i}$ images of marked points lying at $i$-th
boundary contour  ordered according the orientation of that boundary contour.

Fourth, according the definition, the structure of Cardy-Frobenius algebra $H$
includes commutative equipped Frobenius algebra $A$, equipped Frobenius algebra $B$,
homomorphism $\phi: A \to B$ and element $U\in A$. Linear forms
$l_A:A\to\mathbb{C}$ and $l_B:B\to\mathbb{C}$ are those that used in
the definition of equipped Frobenius algebras (see section \ref{definition_c}).

Denote by $K_A\in A$ and $K_B\in B$ Cazimir elements of $A$ and $B$ respectively.
$K_A^g$ is just $g$-th power of the element in algebra $A$.
All those elements for the algebra $H$ were described in section \ref{equipped}.

\begin{theorem}

(i) Let $\Omega$ be an orientable surface ($\varepsilon=1$).
Then we have:
$$\mathcal{H}(\Omega,\{\alpha_p\},\{\beta_q\})=
l_B(\phi(\alpha_1...\alpha_mK_A^g)\beta^1_1\dots\beta^1_{m_1}K_B
\beta^2_1\dots\beta^2_{m_2}\dots K_B\beta^s_1\dots\beta^s_{m_s})$$

(ii) Let $\Omega$ be a non-orientable marked surface.
Then we have:
$$\mathcal{H}(\Omega,\{\alpha_p\},\{\beta_q\})=
l_B(\phi(\alpha_1.....\alpha_mU^{2g})\beta^1_1\dots\beta^1_{m_1}K_B\beta^2_1\dots
\beta^2_{m_2}\dots K_B\beta^s_1... \beta^s_{m_s})$$
\end{theorem}

\proof The theorem follows from \cite{AN}, theorem 4.4. $\Box$

\end{document}